\documentclass[10pt]{article}
\usepackage[utf8]{inputenc}
\usepackage[usenames, dvipsnames]{xcolor}
\usepackage{graphicx}
\usepackage{amsmath,amsfonts,amssymb,textcomp,gensymb}
\usepackage[authoryear,round]{natbib}
\usepackage{bbm}
\usepackage{url}
\usepackage{multirow}
\usepackage{algorithm,algorithmicx,algpseudocode}
\usepackage{fullpage}
\usepackage{stmaryrd}

\newcommand{\Rset}{\mathbb{R}}
\newcommand{\Nset}{\mathbb{N}}
\newcommand{\x}{\mathbf{x}}
\newcommand{\X}{\mathbf{X}}

\newcommand{\Xset}{\mathbb{X}}
\newcommand{\B}{\mathcal{B}}
\newcommand{\class}{\mathcal{C}}
\newcommand{\model}{\mathcal{M}}

\newcommand{\esp}{\mathbb{E}}
\newcommand{\prob}{{\mathbb{P}}}
\newcommand{\cvar}{{\text{CVaR}}}
\newcommand{\one}{\mathbbmss{1}}

\begin{document}

 
\title{Sunflower phenotype optimization under climatic uncertainties using crop models}
\author{Victor Picheny\footnote{INRA, UR875 Mathématiques et Informatique Appliquées Toulouse, F-31326 Castanet-Tolosan, France}\footnote{Corresponding author, victor.picheny@toulouse.inra.fr}, Ronan Tr\'epos$^*$, Bastien Poublan$^*$, 
Pierre Casadebaig\footnote{INRA, UMR AGIR, F-31326 Castanet-Tolosan, France}}

\maketitle

\paragraph{keyword}  Clustering; multi-objective optimization; subset sampling 

\begin{abstract}
Accounting for the annual climatic variability is a well-known issue for simulation-based studies of environmental models. It often requires intensive sampling 
(e.g., averaging the simulation outputs over many climatic series), which hinders many sequential processes, in particular optimization algorithms. We propose here an approach based on a subset selection
of a large basis of climatic series, using an ad-hoc similarity function and clustering. 
A  non-parametric reconstruction technique is introduced to estimate
accurately the distribution of the output of interest using only the subset sampling. 
The proposed strategy is non-intrusive and generic (i.e. transposable to most models with climatic data inputs),
and can be combined to most ``off-the-shelf'' optimization solvers. We apply our approach to sunflower phenotype optimization using the crop model SUNFLO.
The underlying optimization problem is formulated as multi-objective to account for risk-aversion.
Our approach achieves good performances even for limited computational budgets, outperforming significantly more ``naive'' strategies.
\end{abstract}

\section{Introduction}\label{sec:introduction}

Using numerical models of complex dynamic systems has become a central process in sciences.
In agronomy, it is now an essential tool for water resource management, adaptation of anthropic or natural systems to a changing
climatic context or the conception of new production systems.
In particular, in the past two decades crop models have received a growing attention
\citep{boote1996potential,brisson2003overview,brun2006working,bergez.13,brown2014plant,mcnider2014integrated},
as they can be used to help improve the plant performances, either through cultural
practices \citep{grechi.12, wu.12} or model-assisted plant breeding \citep{semenov.13, semenov.14, quilot.12}.

Many times, the objective pursued by model users amounts to solving an optimization problem, that is, find the set of input parameters of the model
that maximize (or minimize) the output of interet (cost, production level, environmental impact, etc.).
Examples of such problems abound with environmental models, including water distribution systems design \citep{tsoukalas2014multiobjective},
agricultural watershed management \citep{cools2011coupling} or
the adaptation of cultural practices to climate change \citep{holzkamper2015assessing}.
In phenotype optimization, \citep[or \textit{ideotype design},][]{martre.15}, plant performance (e.g., yield)
is maximized with respect to its morphological and/or physiological traits.

Within the wide range of potential approaches to solve such optimization problems, \textit{black-box optimization methods}
have proved to be popular in this context \citep{maier2014evolutionary,martre.15,quilot.12}, as they only require
limited expertise in optimization while being quite user-friendly,
as they are in essence non-intrusive (i.e., they only require evaluations of the model at hand).

However, a well-known difficulty, shared by many models users, is to deal with climatic information.
Many agricultural or ecological models require yearly series of day-to-day measures of precipitation, temperature, etc.,
as input variables. This is particularly crucial for agricultural or ecological models, for which
the climate has a preponderant impact on the system.
To avoid drawing conclusions biased by the choice of a particular set (i.e., year) of climatic data, one may either
use scenarii approaches (duplicate the analysis for a small number of distinct climates), or average
the model outputs over a (large) number of climatic datasets. Due to the complex plant-climate interaction,
identifying scenarii may prove to be a very challenging task, and the alternative relies on intensive computation,
which rapidly becomes computationally prohibitive if the analysis is embedded in a loop, even for moderately complex models.

A natural solution is to treat the climate as a random variable, which allows the use of the robust (or noisy) optimization
framework. However, if readily available codes abound for continuous, box-constrained parameters and  deterministic outputs,
solutions become scarce for systems depending on stochastic phenomena.
Besides, the problem formulation becomes more complex, as typically risk-aversion preferences need to be accounted for.

The methodological objective of this paper is two-fold. First, we wish to propose a clear optimization framework
for optimization under climatic uncertainties, and in particular to account for risk-aversion concepts
in a transparent manner.
Second, as both optimization and uncertainty analysis are computationally intensive tasks,
we need to provide an algorithmic solution to solve the problem in reasonable time.
In addition, we wish to remain non-intrusive and generic (i.e. transposable to most models with climatic data inputs).
Finally, in order to facilitate the use of parallel computing, we aim at limiting the
complexity of the algirthm to its minimum.

In this work, we focus on the problem of sunflower ideotype design using the SUNFLO crop model.
SUNFLO is a process-based model which was developped to simulate the grain yield and oil concentration as a function of time,
environment (soil and climate), management practice and genetic diversity \citep{casadebaig.11}.
It allows to assess the performance of sunflower cultivars in agronomic conditions.
A cultivar is represented by a combination of eight genetic coefficients (see Table \ref{tab:sunfloPheno}),
which are the variables to be optimized. The SUNFLO model computes the annual yield $y$ (in tons per hectare)
for a given climatic year.

The rest of this paper is organized as follow: Section \ref{sec:pbdefinition} briefly reviews previous works on phenotype optimization,
describes the SUNFLO model and the multi-objective optimization formulation to solve the problem at hand.
Section \ref{sec:algorithm} is dedicated to the optimization algorithm, which relies on a
subset selection of the available climate data combined with a metaheuristic algorithm.
Finally, Section \ref{sec:results} provides numerical results and compare our approach to
classical solutions.

\section{Problem definition}\label{sec:pbdefinition}
\subsection{Brief review of phenotype optimization}
\citet{martre.15} provide a review of recent developments in this research
domain named \textit{model-assisted crop improvement} or \textit{ideotype design}.
A phenotype is defined as the expression in a particular environmnent of
a specific genotype through its morphology, development, cellular, biochemical
or physiological properties. An ideotype is defined as a combination of morphological
and/or physiological traits
optimizing crop performances to a particular biophysical environment and crop management.
\citet{letort.08} developped an approach to design plant ideotypes maximizing yield,
using numerical optimization methods on coupled genetic and ecophysiological models.
However, as most of the developped crop model do not include genetic-level inputs \citep{hammer2010adapting},
optimization mainly targets the phenotype level.

In the phenotype optimization setting, ideotype design can be formulated
as a problem of optimizing model inputs related to cultivar practices
\citep{grechi.12, wu.12}, or phenotypic parameters \citep{semenov.13, semenov.14, quilot.12}.
Different purposes are targeted such as the
adaptation to climate change \citep{semenov.13, semenov.14}
or the multicriterion assessment of cultivar \citep{quilot.12, qi.10}. In most
of these approaches \citep{letort.08, qi.10, quilot.12}, the study has
been performed on a constant environment, in particular, using a single climatic year.
\citet{quilot.12} stated that further methodological developments are needed in
the  optimization side to reduce computational time in order to be able to consider
multi-environments and large climatic series. In this work, the authors used
the 'Virtual Fruit' model \citep{quilot.05} to
design peach phenotypes for sustainable productions systems. Their aim is to
optimize jointly three model outputs (fruit mass, sweetness and crack density)
in four different scenarii using one climatic data serie in 2009.
They first performed a sensitivity analysis in order to select six phenotypic model 
inputs amongst 60 and use the
multi-objective optimization method NSGA-II \citep{deb.02} in order to solve the problem.

\citet{semenov.13} proposed to evaluate a phenotype by estimating an expected yield using
100 climatic series, by combining the use of the stochastic weather
generator LARS-WG \citep{semenov.10} and the wheat crop model Sirius
\citep{jamieson.98} in
order to design high-yielding ideotypes for a changing climate in the case of
two contrasting situations: Sevilla in Spain and Rothamsted in the United Kingdom.
Inputs were nine cultivar-dependant parameters related to the
photosynthesis, phenology, canopy, drought tolerance and root water uptake.
The optimization problem was solved by using an evolutionary
algorithm with self-adaptation \citep[EA-SA,][]{beyer1995toward}.

\subsection{The SUNFLO model}
In this work, we consider the SUNFLO crop model in order to assess the performance of sunflower cultivars in agronomic conditions.
This model is based on a conceptual framework initially proposed by \citep{monteith.77} and now shared
by a large familly of crop models \citep{Keating2003, brisson2003overview, Stockle2003}. In this framework, the daily crop dry biomass growth rate
is calculated as an ordinary differential equation function of incident photosynthetically
active radiation, light interception efficiency and radiation use efficiency.
Broad scale processes of this framework, the dynamics of leaf area, photosynthesis and biomass
allocation to grains were split into finer processes (e.g leaf expansion and senescence, response functions
to environmental stresses) to reveal genotypic specificity and to allow the emergence of genotype $\times$
environment interactions. Globally, the SUNFLO crop model has about 50 equations and 64 parameters
(43 plant-related traits and 21 environment-related).

In this model, a cultivar is represented by a combination of eight genetic coefficients (see Table \ref{tab:sunfloPheno}).
These coefficients describe various aspects of crop structure or functioning: phenology, plant architecture,
response curve of physiological processes to drought and biomass allocation.
The consequence of genetic modifications can be emulated by changing the values of such parameters.
We consider here the design of sunflower cultivars for a given set of cultural practices and
a specific environment. The overall objective is to find a phenotype that maximizes the yield
for the year to come, without knowing in advance the climate data.
We assume that the coefficients can take continuous values between a lower and an upper bound, which
are determined from a dataset of existing cultivars (see Table \ref{tab:sunfloPheno}). We denote $\x \in \Xset \in \Rset^d$ a
particular phenotype, where $d$ is the number of input variables ($d=8$).

\begin{table}
\caption{\label{tab:sunfloPheno}Phenotypic coefficients and the bounds
used for optimization.}
\centering
\fbox{%
\begin{tabular}{*{4}{l}}

\em Symbol & \em Description  & \em Min & \em Max \\\hline
\multirow{ 2}{*}{TDF1}& Temperature sum from emergence to                 & \multirow{ 2}{*}{765} & \multirow{ 2}{*}{907} \\
~                     & the beginning of flowering ($^{\circ}\mathrm{C}$) & ~   & ~ \\ \hline
\multirow{ 2}{*}{TDM3}& Temperature sum from emergence to                 & \multirow{ 2}{*}{1540} & \multirow{ 2}{*}{1830} \\
~                     & seed physiological maturity ($^{\circ}\mathrm{C}$)& ~   & ~ \\ \hline
\multirow{ 2}{*}{TLN} & \multirow{ 2}{*}{Number of leaves at flowering}   & \multirow{ 2}{*}{22.2}  &  \multirow{ 2}{*}{36.7} \\
~                     &  ~                                                & ~   & ~ \\ \hline
\multirow{ 2}{*}{K}   & Light extinction coefficient                      & \multirow{ 2}{*}{0.780} & \multirow{ 2}{*}{0.950} \\
~                     &  during vegetative growth                         & ~   & ~ \\ \hline
\multirow{ 2}{*}{LLH} & Rank of the largest leave                         & \multirow{ 2}{*}{13.5} & \multirow{ 2}{*}{20.6} \\
~                     & of leaf profile at flowering                        & ~   & ~ \\ \hline
\multirow{ 2}{*}{LLS} & Area of the largest leave of                   &  \multirow{ 2}{*}{334} &  \multirow{ 2}{*}{670} \\
~                     & leaf profile at flowering ($cm^2$)                & ~   & ~ \\ \hline
\multirow{ 2}{*}{LE}  & Threshold for leaf expansion               & \multirow{ 2}{*}{-15.6}  & \multirow{ 2}{*}{-2.31} \\
~                     & response to water stress                             & ~   & ~ \\ \hline
\multirow{ 2}{*}{TR}  & Threshold for stomatal conductance                  & \multirow{ 2}{*}{-14.2} & \multirow{ 2}{*}{-5.81} \\
~                     & response to water stress                & ~   & ~
\end{tabular}}
\end{table}

\begin{figure}
\begin{center}
\includegraphics[width=\textwidth]{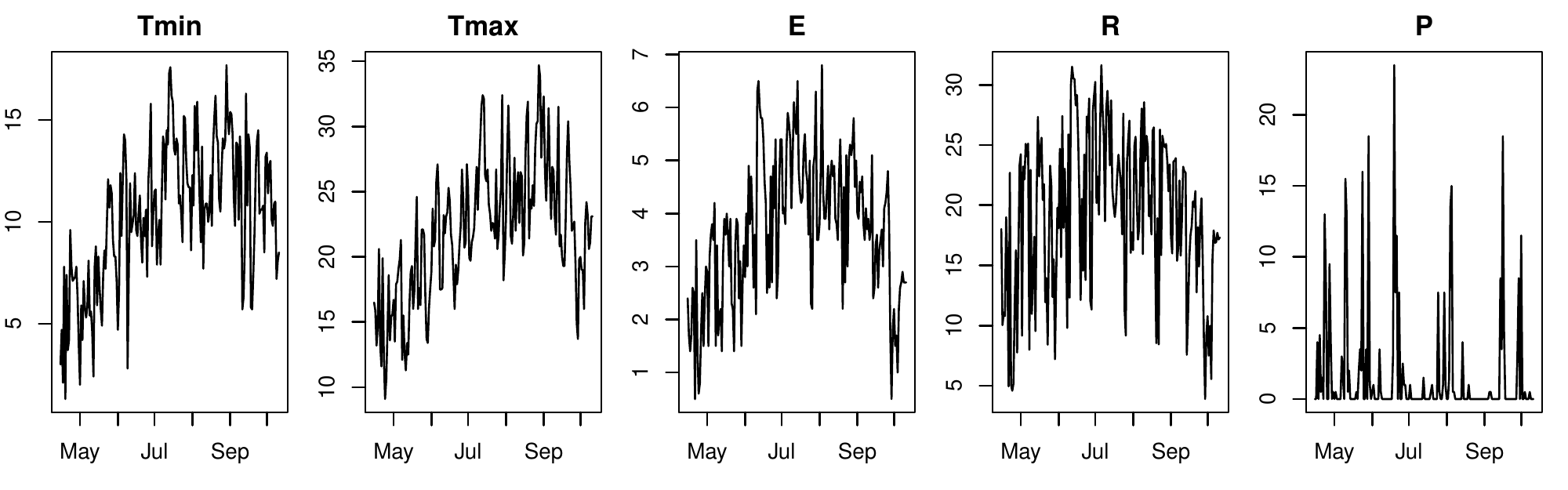}
\caption{\label{fig:climate_example} Dataset of the year 2009, Lusignan.}
\end{center}
\end{figure}

The SUNFLO model computes the annual yield $y$ (in tons per hectare) for a given climatic year.
Hence, it requires as an additional input a climatic serie, which consists of daily measures over a year of five variables:
minimal temperature ($T_{\min}$, \degree Cd), maximal temperature ($T_{\max}$, \degree Cd),
global incident radiation ($R$, $MJ/m^2$), evapotranspiration ($E$, mm, Penman-Monteith) 
and precipitations ($P$, mm)
We note: $c=\{T_{\min}, T_{\max}, R, E, P\}$. Figure \ref{fig:climate_example} provides an example of such data.

\begin{figure}
\begin{center}
\includegraphics[trim=2mm 2mm 10mm 10mm, clip, width=.5\textwidth]{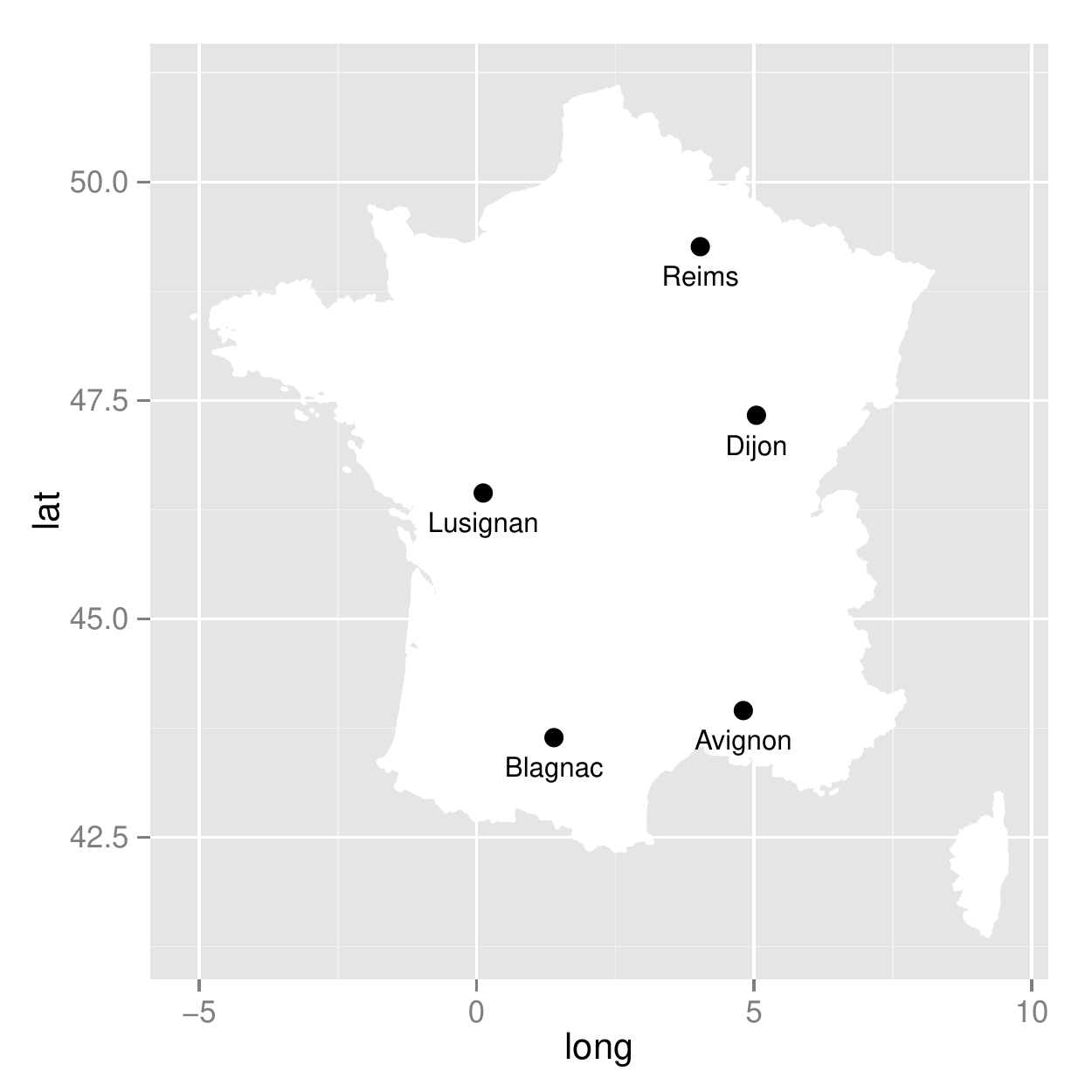}
\caption{\label{fig:climat_france} Location of the five French stations for the
historic climatic data}
\end{center}
\end{figure}

We use historic climatic data
from five french locations Avignon, Blagnac, Dijon, Poitiers and Reims
(see Figure \ref{fig:climat_france}) from 1975 to
2012. The initial data is recorded over 365 days, but we consider only the
cultural year (April to October, 180 days), as the yield computed by the model
only depends on this period.
We denote by $\Omega$ this set of climatic series, and we have $Card(\Omega)=N=190$
and $c \in \Rset^{5 \times 180}$.

To summarize, the yield can be seen as a function of the phenotype and the climatic serie:
\begin{equation*}
 \begin{array}{ccc}
  y: \Xset \times \Omega &\rightarrow &\Rset^+ \\
  \x,c &\longmapsto & y(\x,c).
 \end{array}
\end{equation*}

With a slight abuse of notations, we also define:
\begin{eqnarray*}
 y(\X,c) &:=& \left[y(\x_1,c), \ldots, y(\x_I,c)\right]^T,\\
 y(\x,\mathbf{C}) &:=& \left[y(\x,c_1), \ldots, y(\x,c_J)\right]^T,\\
 y(\X,\mathbf{C}) &:=& \left( y(\x_i,c_j)\right)_{1\leq i \leq I, 1\leq j \leq J},
\end{eqnarray*}
%
that is, the yield function for a set of inputs, either for a set of phenotypes $\X=\{\x_1,\ldots,\x_I\}$ ($I \in \Nset^*$),
a set of climatic series $\mathbf{C}=\{c_1,\ldots,c_J\}$ ($1 \leq J \leq N$), or both.

\subsection{A multi-objective optimization formulation for robust optimization}
The objective is to find a phenotype that maximizes the yield for the year to come, without knowing in advance the climate data.
Let $C$ be the climatic serie of the upcoming year (the upper case denoting a random variable); we consider in the following that $C$ is uniformly distributed
over the discrete set $\Omega$. Since $C$ is random, the yield $y(\x,C)$ is also a random variable (which we denote in the following $Y(\x)$),
which makes its direct maximization with respect to  $\x$ meaningless.

A natural formulation is to maximize the yield expectation:
\begin{equation*}
 \max_{\x \in \Xset} \esp \left[ y(\x,C) \right] = \max_{\x \in \Xset} \esp \left[ Y(\x) \right],
\end{equation*}
with here: $ \esp \left[ Y(\x) \right] = \frac{1}{N} \sum_{i=1}^N y(\x,c_i)$.

However, in general, a farmer also wishes to integrate some prevention against risk in its decision.
Such a problem is often referred to as \textit{robust optimization} in the engineering literature
\citep[see for instance][for a review]{beyer2007robust}.

A popular solution is to replace the expectation by a performance indicator that provides a trade-off
between average performance and risk aversion: typically, the expectation penalized by the variance or
a so-called \textit{utility function}. The drawback of such approaches is that the trade-off must be tuned beforehand
by choosing penalization parameters specific to the method. Choosing the appropriate trade-off may not be straightforward,
and modifying it requires to restart the entire optimization procedure.

We propose here an alternative, which is to consider this problem as multi-objective, by introducing a second criterion to maximize
that accounts for the risk \citep[as in][for instance]{tsoukalas2014multiobjective}. One may choose for instance to maximize a quantile:
\begin{equation*}
 \max_{\x \in \Xset} Q_\alpha \left[ Y(\x) \right],
\end{equation*}
with the usual definition of the quantile: $\prob \left[ Y \leq Q_\alpha(Y) \right] = \alpha$, and $\alpha \in (0, 0.5]$.
Here, it amounts to maximizing the yield for the $(N \times \alpha)$-th worst year.
However, we consider here a close but numerically more stable criterion, called the conditional value-at-risk \citep[CVaR, ][]{rockafellar2000optimization}, defined as:
\begin{equation*}
 \cvar_\alpha \left[ Y(\x) \right] = \esp \left[ Y(\x) | Y(\x) \leq Q_\alpha \left[ Y(\x) \right] \right].
\end{equation*}
$\cvar_\alpha$ is the average yield over the $(N \times \alpha)$-th worst years.

The multi-objective optimization problem is then:
\begin{equation*}
\left\{
\begin{array}{ll}
\max & \esp \left[ Y(\x) \right] \\
\max & \cvar_\alpha \left[ Y(\x) \right] \\
\text{s.t.} & \x \in \Xset.
\end{array}
\right.
\end{equation*}

Such a formulation is relatively classical in robust optimization, although the second objective is often
taken as the variance of the response: $var[Y(\x)]$ \citep[as for instance in][]{chen1999quality,jin2003trade}.
However, considering an expectation-variance trade-off does not make sense here, as a farmer will not want to
minimize the variability of its income (i.e., minimizing the variance) but rather minimize the risk of low income.

\section{Optimization with a representative subset}\label{sec:algorithm}
The two objective functions, $\esp[Y(\x)]$ and $\cvar_\alpha [Y(\x)]$, require running the SUNFLO simulator $N$ times everytime a new phenotype $\x$ is evaluated.
Embedded in an optimization loop, which typically requires thousands to millions
calls to the objective functions, this evaluation step becomes prohibitive.

We propose to address this problem by replacing the large climatic data set $\Omega$ by a small representative set $\Omega_K$. To do so, we first 
choose the set $\Omega_K$ prior to optimization using a clustering algorithm described in Section \ref{sec:classification}. Then, the optimization
algorithm is run using $\Omega_K$. Hence, $\esp[Y(\x)]$ and $\cvar_\alpha [Y(\x)]$ are replaced by their estimates based on $\Omega_K$, which are described
in Section \ref{sec:reconstruction}.

\subsection{Choosing a representative subset of climatic data}\label{sec:classification}

\subsubsection{Principle}
To select our subset, we propose to define a distance (or, conversely, a similarity) between two climatic series,
then choose series \textit{far from each other} using clustering algorithms. 

One can choose to consider only the dataset and define a distance that characterizes differences between the time series. 
However, the drawback of this method is that it is completely model-independent: two climatic series 
can be considered as far from each other but have a similar effect on the model, hence return a similar yield. Inversely, two climatic series can be generally 
close but return different yields because of small critical differences (say, a rainy week at an appropriate moment of the plant growth).

An alternative is to consider a model-based distance: two climatic series would be far from each other only if they return different yields for a given phenotype.
This naturally implies that all the climatic series are run on a (small) phenotype learning set. Therefore, the distance will be very dependent on the choice of 
the set and may result in poor robustness.

Therefore, we propose here to combine both ideas, and define a hybrid distance
that depends on intrinsic differences and on the effect on the model.

\subsubsection{Dissimilarity between time series}
As a climatic serie is defined by five time series of different nature, we need first to define a metric to compare each series separately.
Due to the nature of the data, Euclidian distance can be ruled out, as it makes little sense here. 
Indeed, all the series have important day-to-day variations (corresponding to good or bad weather), and similar events
can be observed from one series to another shifted by one or several days. This is particularly apparent for the precipitation series,
which contain many zeros and several ``peaks'': Euclidian distance would consider two series as far from each other, as long as the 
peaks do not coincide exactly.

A classical tool for time series analysis, sensible in our case, is an algorithm called dynamic time warping 
\citep[DTW, ][]{berndt1994using,aach2001aligning,kadous1999learning}. 
In short, DTW allows two time series that are similar but locally out of phase to align in a non-linear manner, by matching events within a given window.
Note that the DTW algorithm  has a $\mathcal{O}(n^2)$ time complexity, which makes the dissimilarity computation non-trivial. However, this step should be performed only once.
Given two weather series $c_i$ and $c_j$, five distances can be computed, according to the weather variables:
 $d(c_i, c_j)^{Tmin}$, $d(c_i, c_j)^{Tmax}$, $d(c_i, c_j)^{R}$, 
 $d(c_i, c_j)^{E}$ and  $d(c_i, c_j)^{P}$. 

\begin{figure}
\begin{center}
\includegraphics[trim=0mm 5mm 10mm 20mm, clip, width=\textwidth]{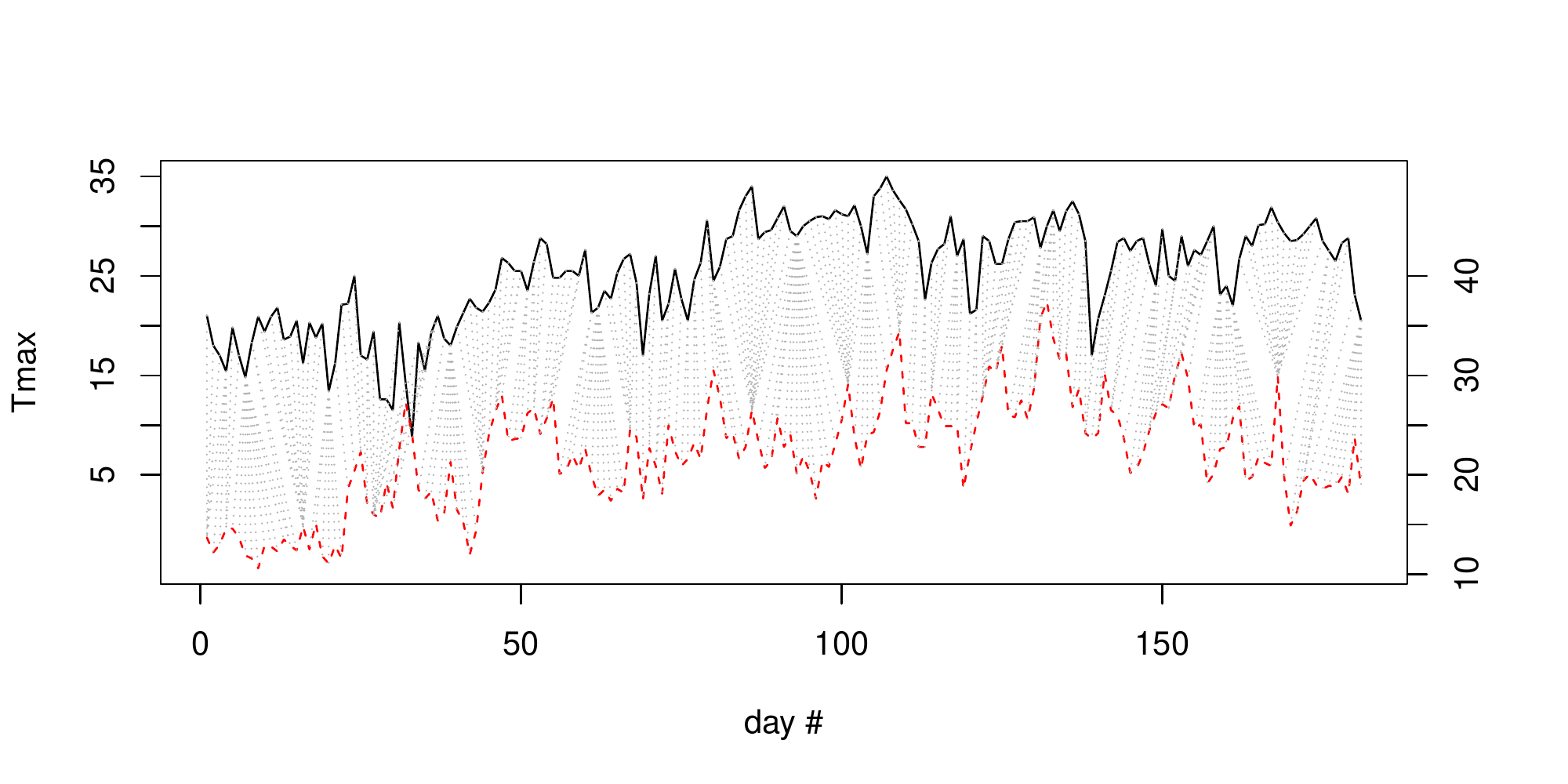}
\caption{\label{fig:dtw} Computing of the dtw distance between two time series
of maximal daily temperature ($Tmax$): Avignon in 1985 (upper curve, left scale) 
and Lusignan in 2012 (bottom curve, right scale).
Dotted line represents the optimal matching of daily temperature 
computed by dtw, for a window size of 7 days.}
\end{center}
\end{figure}

\subsubsection{Model-based dissimilarity}
This dissimilarity measures a difference in the output of the model (the yield). To do so, we choose first a small set of $l$ phenotypes:
$\B = \{\x_1, \ldots, \x_l\}$. 
Typically, $\B$ can be chosen by Latin Hypercube Sampling \citep[LHS,][]{mckay1979comparison} to ``fill'' the search space $\Rset^d$. 
For this basis, the yield is computed for all the climatic series: $y(\B, \Omega) \in \Rset^{l \times N}$.
Then, the model-based distance is simply the Euclidian distance:
\begin{equation*}
 d(c_i, c_j)^{\model} = \sqrt{ \frac{1}{l} \sum_{k=1}^l \left( y(\x_k, c_i) - y(\x_k, c_j) \right)^2}
\end{equation*}

\subsubsection{Combining dissimilarities}\label{sec:finaldistance}
We want here to combine the six dissimilarities (one for each time series and the model-based one) into a single one, 
with equal weight to each variable. We propose to do so by normalizing the dissimilarities before summing them with uniform weights.
As the variables are of different nature, the dissimilarities distributions are likely to be very different (uniform, heavy tailed, etc.), 
hence artificially weight the variables even if they are rescaled similarly.

Here, we follow a normalization procedure proposed in \citet{olteanu2015line} called ``cosine preprocessing'', which works as follow: 
Let $\mathbf{D}$ be a $N \times N$ matrix of dissimilarities (with values $d_{ij} = d(\x_i,\x_j)$, $d_{ij}=d_{jj}$ and $d_{ii}=0$).
We first compute a corresponding similarity matrix $\mathbf{S}$, with values:
\begin{equation*}
 s_{ij} = -\frac{1}{2} \left[ d_{ij} - \frac{1}{N} \sum_{k=1}^N \left( d_{ik} + d_{kj} \right) + \sum_{k=1}^N \sum_{k'=1}^N d_{kk'} \right].
\end{equation*}
Then, we normalize $\mathbf{S}$ with:
\begin{equation*}
 \bar s_{ij} = \frac{s_{ij}}{\sqrt{s_{ii} + s_{jj} }},
\end{equation*}
and the normalized dissimilarity matrix $\mathbf{\bar D}$ has elements defined as:
\begin{equation*}
 \bar d_{ij} = \bar s_{ii} + \bar s_{jj} - 2\bar s_{ij} = 2 - 2\bar s_{ij}.
\end{equation*}

Now, we use a convex combination of the six normalized dissimilarities:
\begin{equation}
 \delta_{ij} = \alpha_{T_{\min}} \bar d_{ij}^{T_{\min}} + \alpha_{T_{\max}} \bar d_{ij}^{T_{\max}} + \alpha_{P} \bar d_{ij}^{P} + \alpha_{E} \bar d_{ij}^{E} + \alpha_{S} \bar d_{ij}^{S} + \alpha_{\model} \bar d_{ij}^{\model},
\end{equation}
with $\alpha_{T_{\min}} + \ldots + \alpha_{\model} = 1$. In the following, we use $\alpha_{\model} = 1/2$ and the other weights equal to $1/10$.

\subsubsection{Choosing a representative subset using classification}
Once the matrix of dissimilarities $\boldsymbol{\Delta}$ is computed, most unsupervised clustering algorithms can be used to split the set 
of climatic series $\Omega$ into subsets.
However, a difficulty here is that the centroids of the clusters cannot be computed. Hence, we use a variation of the k-means algorithm that only requires \textit{dissimilarities}
to the centroids. We follow the approach described in \citet{olteanu2015line}; the corresponding pseudo-code is given in Algorithm \ref{al:kmeans}.

The algorithm divides the set $\Omega$ into $K$ classes $\class^{(1)}, \ldots, \class^K$, not necessarily of equal sizes. 
A class $\class^k$ contains $N^k$ elements $\{c_1^k, \ldots, c_K^k\}$. 
Any element $c \in \Omega$ is uniquely attributed to one class and we have: $\sum_{k=1}^K N^k = N$.
For each class $k$, a representative element $\omega^k$ is chosen, which we use to define the representative set:
$\Omega_K = \{\omega^1, \ldots, \omega^K \}$.

\subsection{Non-parametric reconstruction of distributions}\label{sec:reconstruction}
\begin{figure}
\begin{center}
\includegraphics[natwidth=10in,natheight=7in,scale=0.8]{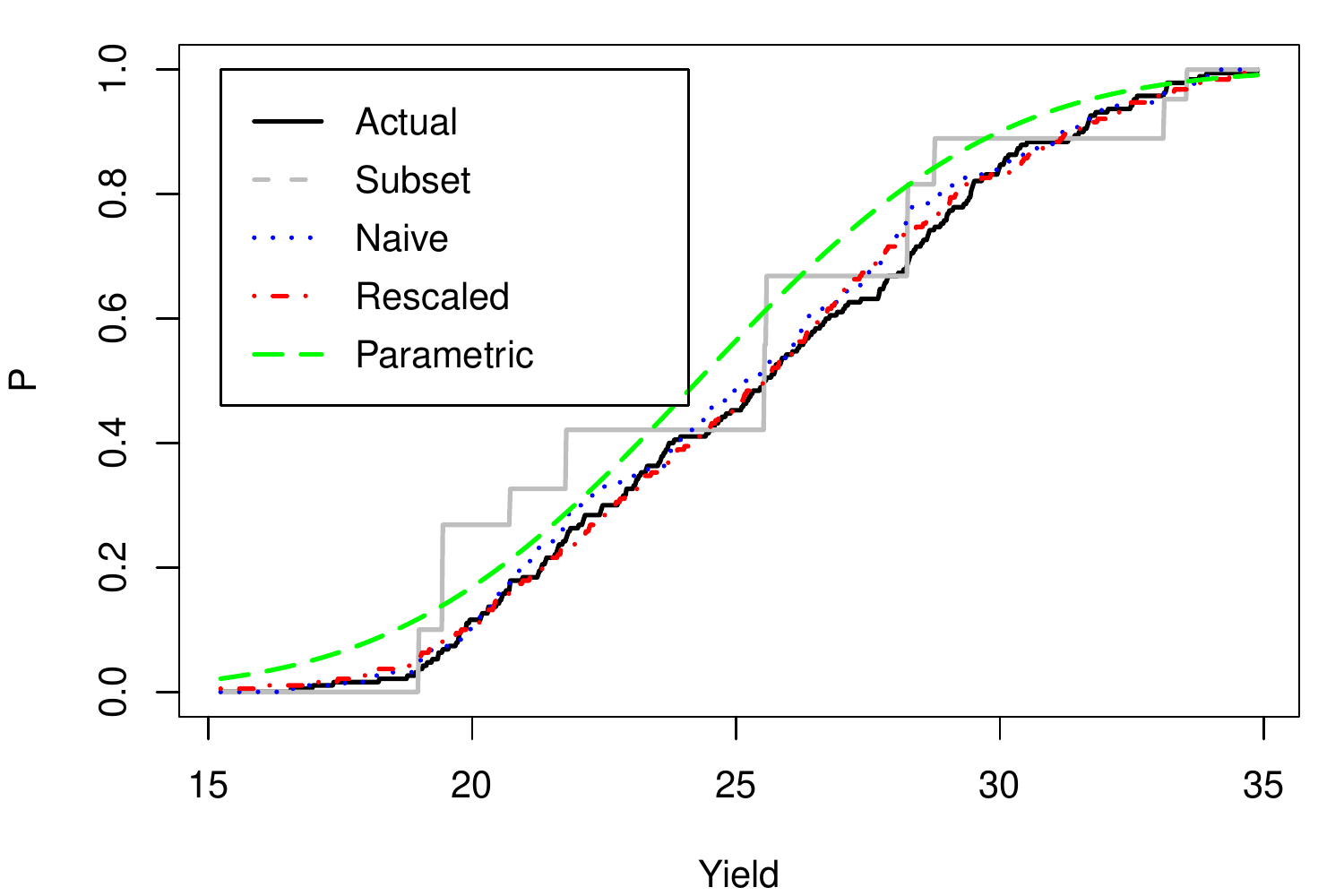}
  \caption{\label{fig:CDF} Actual and estimated distributions (CDF) of the yield of a given phenotype.}
\end{center}
\end{figure}

The objective here is to obtain accurate estimations of the objective functions $\esp[Y(\x)]$ and $\cvar_\alpha [Y(\x)]$ based on 
the yield computed for a new phenotype and the representative set: $y(\x, \Omega_K)$. Since this set is small,
computing directly the objective functions would lead to large errors, in particular for $\cvar_\alpha [Y(\x)]$, 
that requires an accurate representation of the tail distribution (see Figure \ref{fig:CDF}).
A natural alternative is to fit a parametric distribution the small data set, and infer the objectives on the distribution.
However, the form of the empirical distribution (Figure \ref{fig:CDF}) does not readily call for a given parametric model,
and misspecifying the distribution shape may result with large bias.

Hence, we propose to reconstruct the distribution using a non-parametric method, by re-using the 
data computed for the classification step, that is, the yield computed for the phenotype learning basis 
and all the climatic series ($y(\B, \Omega)$).

The general idea is to consider a mixture model for the yield (each component corresponding to a class $\class^k$):
\begin{equation*}
 f_{Y(\x)}(y) = \sum_{k=1}^K \frac{N^k}{N} f_{Y^k(\x)}(y), \qquad y\in\Rset,
\end{equation*}
$f$ standing for the probability density function (PDF), and $Y^k(\x)$ being yield within the class $k$.

We decompose further $Y^k(\x)$ as the sum of the value at the representative element and a residual:
\begin{equation*}
 Y^k(\x) = y(\x, \omega^k) + \varepsilon^k(\x).
\end{equation*}

The intra-class distribution is then characterized by the residuals $\varepsilon^k(\x)$, 
which determine the form, spread (or amplitude), and bias (i.e., difference between the average value and the 
value of the representative element). All these elements vary from one class to another, which advocates the 
use of non-parametric approaches.
%

\begin{figure}
\begin{center}
\includegraphics[width=\textwidth]{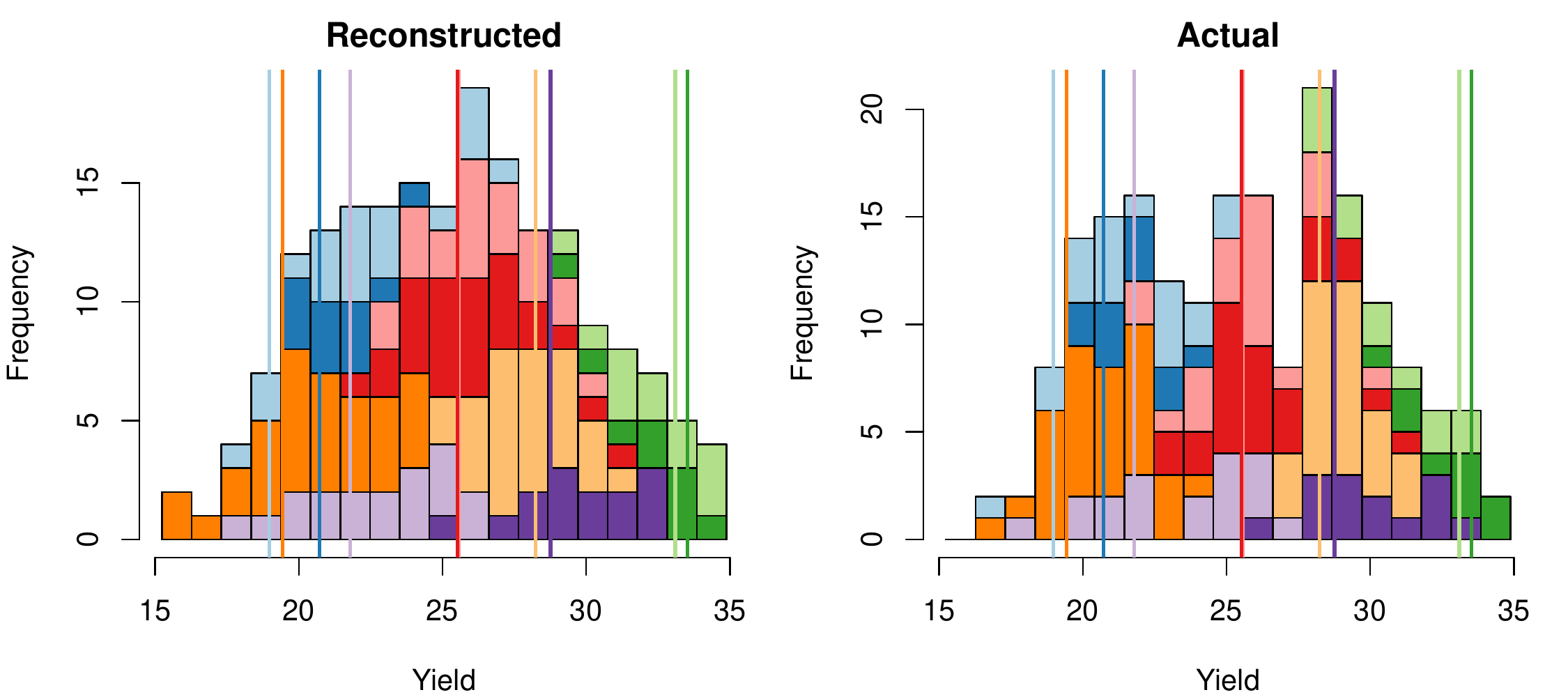}
\caption{\label{fig:histogram} Estimated yield distribution of a given phenotype. The colors show how the reconstruction works: 
  each color corresponds to a class $k$, and the vertical bars to the representative element of the class.}
\end{center}
\end{figure}

\paragraph{Method 1 (naive)}
From $y(\B, \Omega)$, we first compute the residuals $\varepsilon_j^k(\x_i)=y(\x_i, c_j^k) - y(\x_i, \omega^k)$ ($1 \leq i \leq l$; $1 \leq j \leq N^k$; $1 \leq k \leq K$).
Then, we average the residuals over the phenotypes of $\B$:
\begin{equation*}
\hat{\boldsymbol{\varepsilon}}^k = \left[ \hat\varepsilon_1^k, \ldots, \hat\varepsilon_{N^k}^k \right] \text{, with } \hat\varepsilon_j^k = \frac{1}{l} \sum_{i=1}^{l} \varepsilon_j^k(\x_i).
\end{equation*}
The intra-class yield variety is re-created by adding the average residual vector to the yield computed for the representative value:
\begin{equation*}
 {\hat{Y}}^k(\x) = y(\x, \omega^k) + \hat{\varepsilon}_i^k,
\end{equation*}
with $i$ uniformly taken from $\llbracket 1, {N^k} \rrbracket$.
Thus, each component of the mixture has a fixed distribution (i.e. independent of $\x$), shifted according to its representative value, and the mixture 
shape and spread varies according to the distribution of the representative values (see Figure \ref{fig:histogram} for an illustration).

However, in practice, the values of the residuals can vary substantially from one phenotype to another, 
and averaging them over $\B$ tends to destroy the shape information. To address this issue, we proposed the following 
modification:
%

\paragraph{Method 2 (rescaled)} 
We introduce first the weighted variance of the yield over the representative set:
\begin{eqnarray*}
 \sigma_K^2(\x) = \frac{1}{N} \sum_{k=1}^{K}{N^k \left(y(\x, \omega_k) - \frac{1}{N} \sum_{j=1}^{K} N^j y(\x, \omega_j) \right)^2}.
\end{eqnarray*}
Note that for a new phenotype $\x$, the only data available is indeed $y(\x, \omega_j)$, so few alternatives are possible.
We then define averages of \textit{normalized} residuals:
\begin{equation*}
\bar{\boldsymbol{\varepsilon}}^k = \left[ \bar\varepsilon_1^k, \ldots, \bar\varepsilon_{N^k}^k) \right] \text{, with } \bar\varepsilon_j^k = \frac{1}{l} \sum_{i=1}^{l} \frac{\varepsilon_j^k(\x_i)}{\sigma_K(\x_i)}.
\end{equation*}
and the yield vector is constructed with:
\begin{equation*}
 {\hat{Y}}^k(\x) = y(\x, \omega^k) + \sigma_K(\x) \times \bar{{\varepsilon}}_i^k,
\end{equation*}
with $i$ uniformly taken from $\llbracket 1, {N^k} \rrbracket$.

Figures \ref{fig:histogram} and \ref{fig:CDF} illustrate the reconstruction technique for a given (randomly chosen) phenotype.
On Figure \ref{fig:histogram}, we see how the estimated distribution is built using the residuals corresponding to each class.
We can see that the range and shape of the residuals vary considerably from one class to another. Also, their distribution around 
the representative element differs: as the residuals do not have a zero mean, the value of the representative element is not 
necessarily central for each class. Comparing the reconstructed (Figure \ref{fig:histogram}, top) and actual (bottom) distributions,
we see that the mixture is globally the same on both graphs.

Figure \ref{fig:CDF} shows the cumulative distribution function (CDF) of the actual yield and of three estimations: using the two 
methods described above and a simple parametric method, which consists in assuming a Gaussian distribution of the yield. 
The empirical CDF corresponding to the subset values only is also depicted, with unequal steps to account for the different 
number of elements in each class.

We first notice that the subset data only is obviously insufficient to evaluate accurately the mean or the CVaR. Then, we see that
the actual distribution does not seem to belong to a known distribution, and using a normal distribution introduces a large bias.
Inversely, using a non-parametric reconstruction allows us to match the shape of the actual distribution. The difference between
the two methods is small for this example, yet the second approach is slightly better almost everywhere.

In our study, we found that this second method provided a satisfying trade-off between robustness, simplicity and accuracy. 
Yet, many refinements would be possible at this point, for instance by introducing intra-class rescaling (different normalization for each class), 
bias correction, or using the distance from the phenotype $\x$ to the basis $\B$.

\subsection{Optimization and reconstruction update}\label{sec:improving}
Finally, the multi-objective optimization problem solved is:
\begin{equation*}
\left\{
\begin{array}{ll}
\max & \esp \left[ \hat{Y}(\x) \right] \\ 
\max & \cvar_\alpha \left[ \hat{Y}(\x) \right] \\
\text{s.t.} & \x \in \Xset,
\end{array}
\right.
\end{equation*}
with $\hat{Y}(\x)$ a mixture of $\hat{Y}^1(\x), \ldots, \hat{Y}^K(\x)$.

One may note that $\esp \left[ \hat{Y}(\x) \right]$ and $\cvar_\alpha \left[ \hat{Y}(\x) \right]$ serve as estimates of $\esp \left[ Y(\x) \right]$ and 
$\cvar_\alpha \left[ Y(\x) \right]$, respectively. These estimates are based on the phenotype basis $\B$, which is sampled uniformly over $\Xset$
to offer a general representation of the phenotype space. This feature is important at the beginning of the optimization to ensure that the optimizer 
does not get trapped into poorly represented regions. However, as the optimizer converges towards the solution, the search space becomes more narrow, 
and a substantial gain in performance can be achieved by modifying the estimates so that they are more accurate in the optimal region.

In theory, it is possible to re-run the entire clustering procedure after a couple of optimization iterations, by adding new phenotypes to the learning set. 
However, such strategy is likely to increase greatly the computational burden. We propose instead to modify only the reconstruction step, for which only
very few additional calculations are required.

Indeed, the reconstructed yield ditributions use the phenotype learning basis $\B$ and their associated values $y(\B, \Omega)$.
By replacing the initial $\B$ with $\B'$ formed by phenotypes chosen inside the optimal region, we obtain yield values $y(\B',\Omega)$
that are more likely to represent the actual distribution within this region. Such ``specialization'' may be to the detriment 
of the global accuracy of the estimates, but this is not critical as the optimizer concentrates on a narrow region.

Including a new phenotype $\x'$ into the basis $\B$ requires running the SUNFLO simulator $N$ times to obtain $y(\x', \Omega)$.
Therefore, an efficient trade-off must be found between pursuing the optimization and improving the estimates. Also, it may be 
beneficial to discard phenotypes in $\B$ that are far from the optimal region. In summary, we need to:
a) decide when to add phenotypes to the basis and b) when to discard them and c) choose which to add / discard.

A simple strategy is to perform only two steps: first, run the optimization with the initial basis $\B$. Then, select $l$ new phenotypes from the 
obtained Pareto set and replace the entire basis $\B$ after running the $N \times l$ simulations. Finally, restart the optimization with the new estimates.
We have found (Section \ref{sec:results}) that this two-step strategy was sufficient on our problem, while relatively easy to implement.

\subsection{Optimization procedure overview}\label{sec:overview}
To summarize this section, Algorithm \ref{al:general} describes the complete optimization procedure, 
including the initial clustering  and the two-step strategy.
Each step relies on the call to a metaheuristic algorithm such as NSGA-II 
\citep{deb.02} or MOPSO-CD \citep{raquel.GECCO05}. Hence, two-step MOPSO-CD
stands for the tow-step algorithm using the MOPSO-CD metaheuristic.

\begin{algorithm}[H]
 \caption{Two-step optimization algorithm}
 \label{al:general}
 \begin{algorithmic}
  \State \textbf{Initialization}
  \State Choose phenotype database $\B$, and compute yield matrix $y(\B, \Omega)$
  \State Compute matrix of dissimilarity $\boldsymbol{\Delta}$
  \State Run clustering algorithm to obtain the classes $\class^{(1)}, \ldots, \class^K$ and the representative set $\Omega_K$
  \State Get residuals from $y(\B, \Omega)$
  \State
  \State \textbf{Optimization: run 1 }
  \State Choose population size $q$ and number of iterations $T$
  \For {t= 1, \ldots, T}
  \State Select new phenotypes $\{\x_{new}^1, \ldots, \x_{new}^q\}$ according the metaheuristic.
  \State Calculate yield for the representative set for each new phenotype $y(\x_{new}^i, \Omega_K)$
  \State Reconstruct $\hat{Y}(\x_{new}^i)$ and evaluate $\esp \left[ \hat{Y}(\x_{new}^i) \right]$ and $\cvar_\alpha \left[ \hat{Y}(\x_{new}^i) \right]$
  \State Post process $\{\x_{new}^1, \ldots, \x_{new}^q\}$ according the metaheuristic.
  \EndFor
  \State Get Pareto-optimal solutions $\X^*$
  \State
  \State \textbf{Optimization: run 2 }
  \State Replace $\B$ by $\X^*$, compute yield matrix $y(\B, \Omega)$
  \State Get the new residuals from $y(\B, \Omega)$
  \For {t= 1, \ldots, T}
  \State Select new phenotypes $\{\x_{new}^1, \ldots, \x_{new}^q\}$ according the metaheuristic.
  \State Calculate yield for the representative set for each new phenotype $y(\x_{new}^i, \Omega_K)$
  \State Reconstruct $\hat{Y}(\x_{new}^i)$ and evaluate $\esp \left[ \hat{Y}(\x_{new}^i) \right]$ and $\cvar_\alpha \left[ \hat{Y}(\x_{new}^i) \right]$
  \State Post process $\{\x_{new}^1, \ldots, \x_{new}^q\}$ according the metaheuristic.
  \EndFor
  \State Get Pareto-optimal solutions $\X^*$
 \end{algorithmic}
\end{algorithm}

\section{Experimental setup}\label{sec:experiments}
\subsection{Climate subset selection}
In this experiment, we used the \texttt{R} package \texttt{dtw} \citep{giorgino2009computing} 
to compute all the distances between climatic series.
Note that the window size (that is, the maximum shift allowed) is a critical parameter of the method; 
we use here expert knowledge to choose it. For the precipitation, a window of $\pm 3$ days is used; for the other
variables, a window of $\pm 7$ days is chosen. 
The phenotype basis $\B$ is chosen as a 10-point LHS; hence, for this step the method required $1,900$ calls to the SUNFLO model.

Once the dissimilarity matrix $\boldsymbol{\Delta}$ is computed, the clustering algorithm (see Appendix \ref{sec:appendix})
is run. Since this algorithm amounts to a gradient descent, it provides a local optimum only, so we need to restart it several times
(by changing the initial values $\boldsymbol{\beta}_0$) to ensure that a good optimum is found. We found in practice that 
$500$ iterations and $10$ restarts were sufficient to achieve a good robustness.

This algorithm does not choose automatically the number of classes $K$. We found empirically that $K=10$
provided a satisfying trade-off between the representation capability of the subset and the computational cost
during the optimization loop.

\subsection{Optimization}
To solve the multi-objective optimization problem, we chose to use the 
MOPSO-CD metaheuristic \citep[Multi-Objective Particle Swarm Optimization with Crowding Distance][]{raquel.GECCO05}.
MOPSO-CD is a stochastic population-based algorithm inspired by the social behavior of bird flocking. 
In short, the algorithm maintains over $T$ generations a population $P$ of individuals (candidate solutions).
At each generation, each candidate is moved through the search space according to an individual direction (local improvement),
a global direction (towards the best candidates of the population) and a crowding distance. 
This distance is used in order to build a set of solution that fills uniformly the Pareto front.

In the following experiments we used the \texttt{R} package \texttt{dtw} \citep{naval2013mopsocd}.
The two main parameters of MOPSO-CD are the population size and number of generations (their product 
being equal to the number of function evaluatuions).

In order to assess the validity of our approach, we have conducted and empirical
comparison to simpler approaches: random search and a ``naive'' optimizer, both using 
the full set of climatic series. In addition, we have conducted an intensive experiment
to obtain an accurate representation of the actual Pareto set.

The intensive experiment consists in running two multi-objective algorithms 
(NSGA-II and MOPSO-CD) with a very large budget (number of calls to the simulator function)
using the full set of climatic series. The two obtained Pareto fronts are merged to a single one,
which we consider as ``exact'' in the following. We set the number of iterations to 300 and the population size 
to 200, hence computing the exact Pareto front requires 
$2  \times 200 \times 300 \times 190 =22,800,000$ calls to SUNFLO.

Random search, or LHS search, is performed using a latin hypercube sampling approach to 
fill the search space $\Xset$. 
The naive optimization is performed using the original MOPSO-CD algorithm. 
Each sampled point is evaluated using the entire 
set of climatic series ($N=190$) to estimate the expected yield and CVaR. 

We compare the different approaches based on an equal number of calls to SUNFLO (that is, 
we do not consider the time costs related to each approach). We considered four budgets:
large ($380,000$), medium ($95,000$), small ($23,750$) and very small ($11,400$).

For the naive and two-step approaches, we need to define
the number of iterations and the population size. We set the number of iterations  to approximately
five times the popuplation size, except for the very small budget where the population size would be to small. 
For the two-step algorithm, each evaluation of the expectation and CVaR
requires $10$ SUNFLO runs, which allows a larger population and number of iterations 
than the naive approach, but it is also necessary to compute two times $y(\B, \Omega)$ 
(the simulations of yields over all climatic series for the phenotype basis), 
which has a $10 \times 190$ cost.

The different setups are given in Table \ref{tab:expe}. Note that 
the budgets are only approximately equal (due to rounding issues).
Nevertheless the budgets for the two-step approach are always equal or smaller than the naive one.

Since these three optimization approaches are stochastic, each experiment is
replicated 10 times, to assess the robustness of the results.

The time cost of one call to the SUNFLO model is low ($\simeq 0.1$ sec),
which makes it possible to perform such an extensive experiment.
However, to limit the computational costs, these experiments are performed with
either a symmetric multiprocessing (SMP) solution based on 30 cores or a 
message passing interface (MPI) implementation based on 40 cores, depending 
on memory requirements of experiments, which makes time costs comparisons meaningless. 

\begin{table}
\caption{\label{tab:expe}Experiments performed for the two-step MOPSO-CD algorithm evaluation.}
\centering
\begin{tabular}{ |c|cccc| }
\hline
Optimization                & \multirow{2}{*}{Budget} & Nb of                             & Pop                    & Real nb                                      \\
experiment                  & ~                       & iterations                        & size                   & of simulations                    \\ \hline 
\multirow{2}{*}{Intensive}  & Very                    & \multirow{2}{*}{$300 (\times 2)$} & \multirow{2}{*}{$200$} & \multirow{2}{*}{$\simeq2 \times 10^7$}  \\
~                           & large                   & ~                                 & ~                      & ~                                 \\ \hline
~                           & very small              & -                                 & 60                     & 11,400                             \\ \cline{2-5}
Random                      & small                   & -                                 & 125                    & 23,750                             \\ \cline{2-5}
(or LHS)                    & medium                  & -                                 & 500                    & 95,000                             \\ \cline{2-5}
~                           & large                   & -                                 & 2,000                  & 380,000                            \\ \hline
~                           & very small              & 12                                & 5                      & 12,350                             \\ \cline{2-5}
Naive                       & small                   & 25                                & 5                      & 24,700                             \\ \cline{2-5}
MOPSO-CD                    & medium                  & 50                                & 10                     & 96,900                             \\ \cline{2-5}
~                           & large                   & 100                               & 20                     & 383,000                            \\ \hline
~                           & very small              & $42 (\times 2)$                   & 9                      & 11,540                             \\ \cline{2-5}
Two-step                    & small                   & $71 (\times 2)$                   & 14                     & 23,960                             \\ \cline{2-5}
MOPSO-CD                    & medium                  & $152 (\times 2)$                  & 30                     & 95,600                             \\ \cline{2-5}
~                           & large                   & $308 (\times 2)$                  & 61                     & 380,780                            \\ \hline
\end{tabular}
\end{table}

The SUNFLO model has been implemented on the VLE software \citep{quesnel.09}
in the RECORD project which is  dedicated to agorecosystems
study \citep{bergez.13}. VLE is a multi-modeling and simulation platform coded
in C++ that provides both a shared memory and a MPI based parallelisation for
the simulation of multiple input combinations. A native port $rvle$ to the
sofware $R$ is available in order to call simulations from
this statistical tool. The other $R$ packages used are \texttt{fExtremes}
(computation of CVaR statistic), \texttt{lhs} (optimized LHS generation), 
\texttt{emoa} (dedicated tools for multiobjective problems) and \texttt{mco} (NSGA-II implementation). 
Finally, we are grateful to the genotoul bioinformatics platform Toulouse Midi-Pyrenees for providing help
and/or computing and/or storage resources.
\section{Results and discussion}\label{sec:results}
\subsection{Climate subset selection}
We analyze first the classification obtained with our approach.
As the classification is based on non-trivial distances, it is difficult
to characterize each class with integrated quantities (e.g. rainy / hot years, etc.).
We provide in the following three tools for this analysis.

\begin{figure}
\begin{minipage}[t]{0.05\textwidth}
\begin{minipage}[t]{\textwidth}
\vspace{-3cm}
(a)
\end{minipage}
\begin{minipage}[t]{\textwidth}
\vspace{+0.1cm}
(b)
\end{minipage}
\begin{minipage}[t]{\textwidth}
\vspace{+4cm}
(c)
\end{minipage}
\end{minipage}
\begin{minipage}[t]{0.9\textwidth}
\begin{minipage}[t]{\textwidth}
\includegraphics[width=\textwidth]{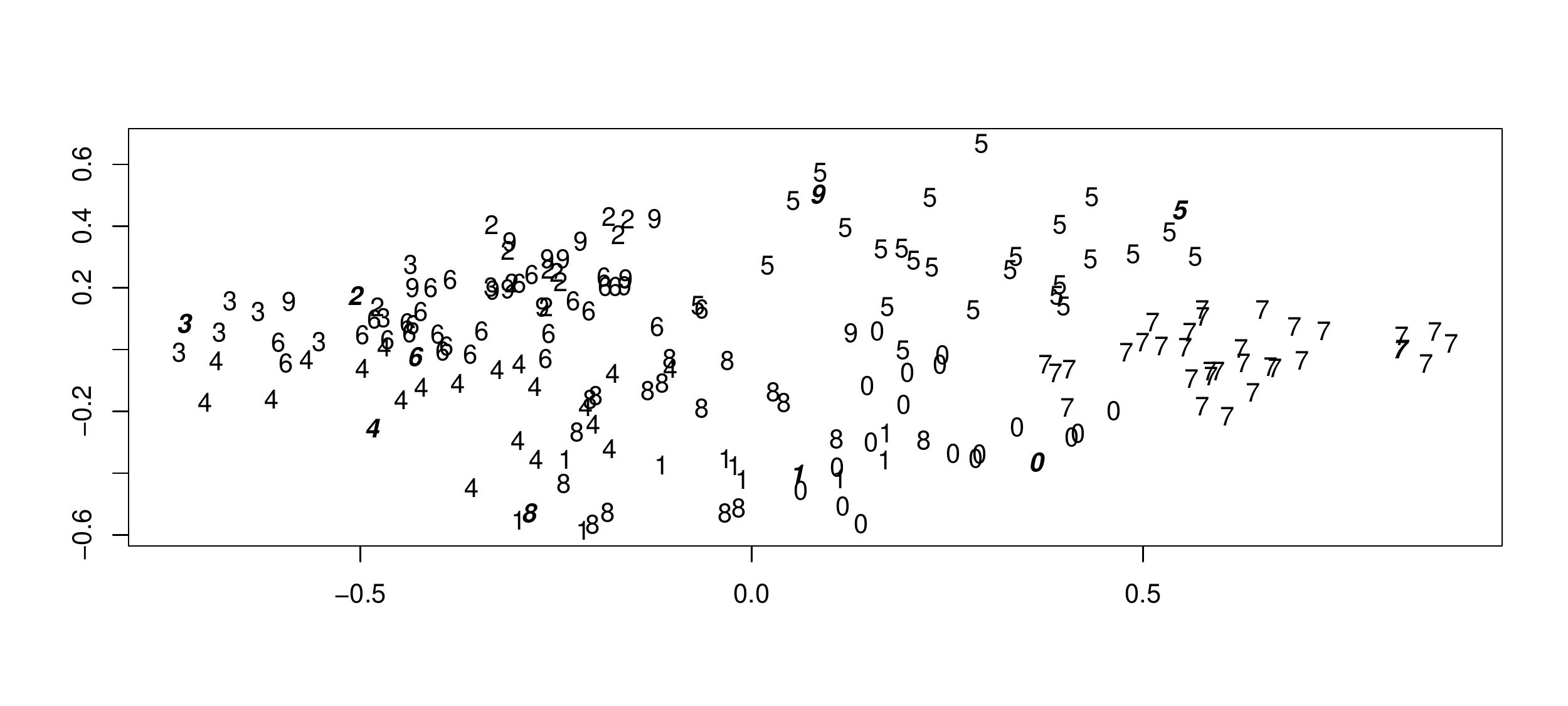}
\end{minipage}
\begin{minipage}[t]{\textwidth}
\includegraphics[trim=0.0cm 0.3cm 1.5cm 3cm, width=\textwidth]{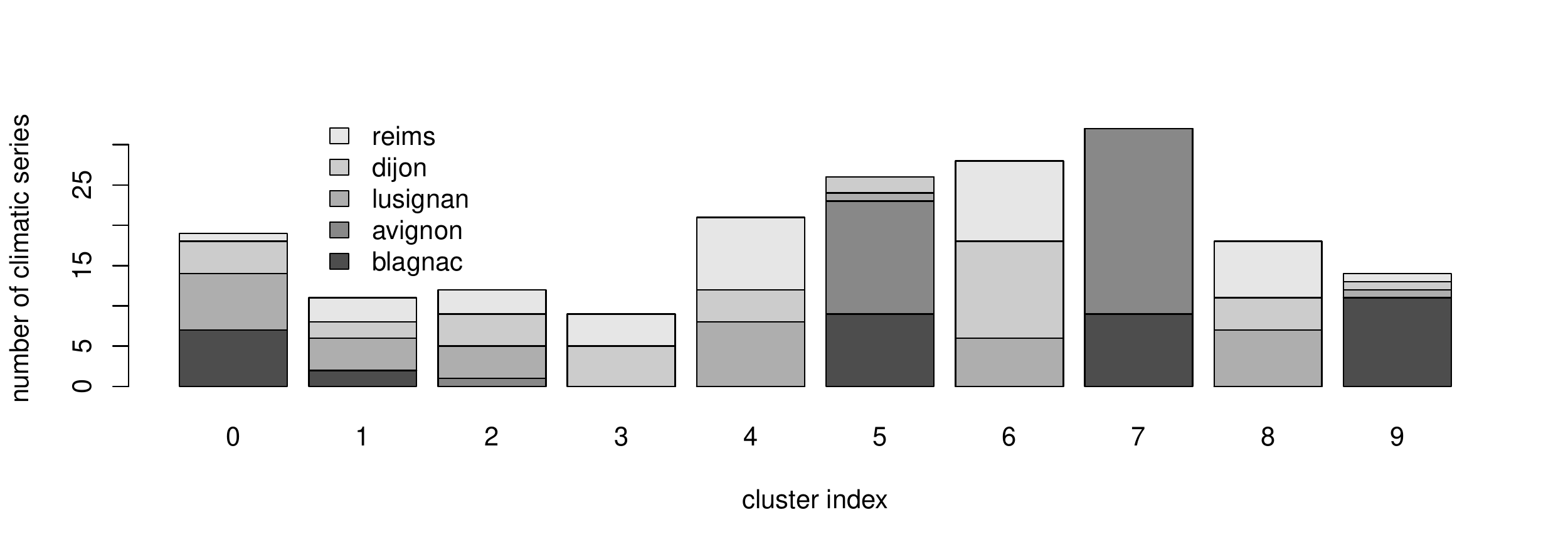}
\end{minipage}
\begin{minipage}[t]{\textwidth}
\includegraphics[width=\textwidth]{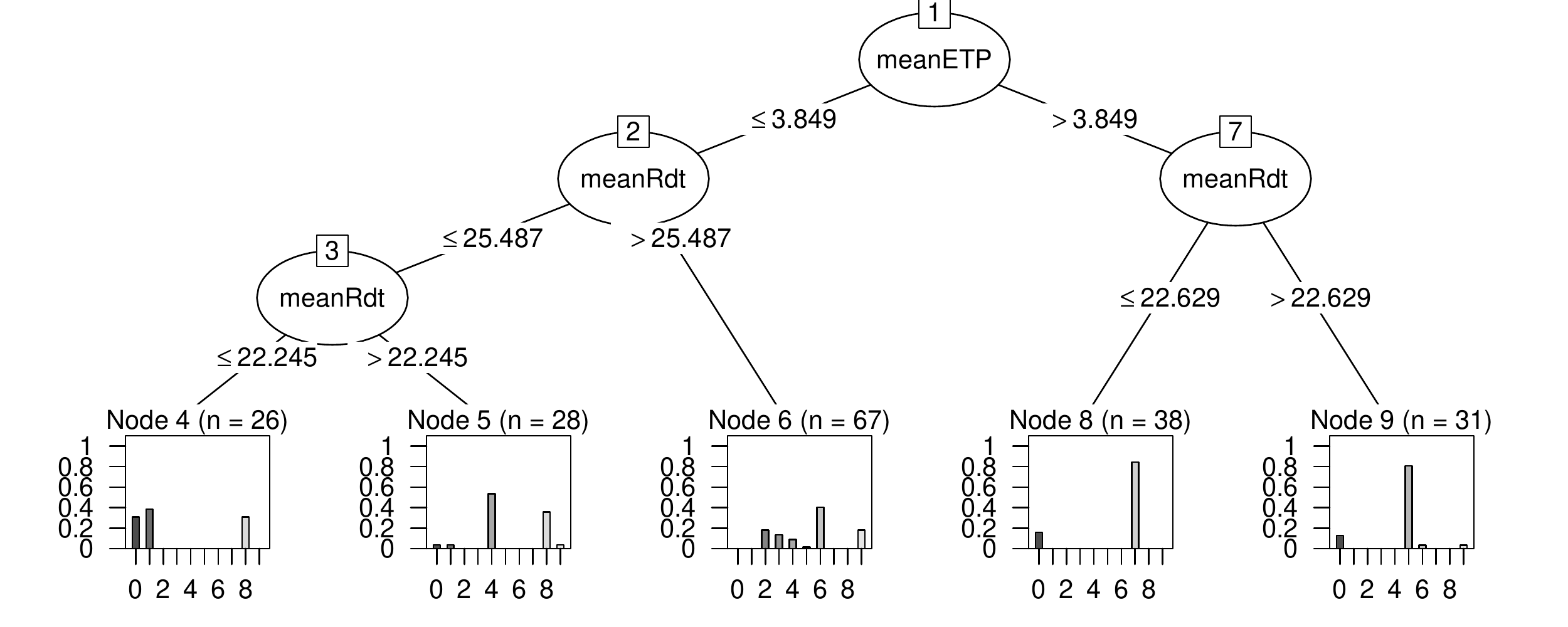}
\end{minipage}
\end{minipage}
\caption{\label{fig:classif_caracterisation} (a) Clusters and individuals (the
190 time series) are plotted in a 2D projection using Classical Multidimensional Scaling.
Each digit represents a weather time series which value corresponds to its
cluster. Climate series of representative set $\Omega_K$ are plotted in bold italic.
(b) Number of climatic series by cluster splitted by
localization in France.
(c) A decision tree to explain clusters using, for each climatic series,
  the temporal mean values of climatic variables and
 the mean yield simulated on the 10 phenotypes in $\B$.}
\end{figure}

We first plot a 2D projection of the climatic series based on the matrix
of distances $\boldsymbol{\Delta}$ computed as in Section \ref{sec:finaldistance}.
To do this, we use the \texttt{R} package \texttt{cmdscale} (Classical Multidimensional Scaling)
(Figure \ref{fig:classif_caracterisation}-a).
Such a representation allows us to see whereas the classes are well-separated, if there are outliers, etc.

In Figure \ref{fig:classif_caracterisation}-b, the number of climatic series,
grouped by their localization is given for each cluster. Finally, a decision tree
has been learnt (with the \texttt{R} package \texttt{C50})
using the cluster index of climatic series as the variable to explain (Figure
\ref{fig:classif_caracterisation}-c). 
We highlight here that this tree is solely for interpretation purpose and
is not linked to the proposed classification strategy.
We used temporal mean aggregation of climatic
 variables $\{T_{\min}, T_{\max}, R, E, P\}$ and the mean yield simulated on the
 10 phenotypes in $\B$ to build the decision tree.

Based on these three representations, one can conclude that some clusters
correspond more or less to wheater types from the South of France (Avignon, Blagnac : 0, 5, 7, 9)
rather warm (5, 7) or not (0, 9) and leading to high yields (5, 9)
or not (0, 7). The three clusters 0, 5 and 7 seem indeed the most easy to characterize
(Figure \ref{fig:classif_caracterisation}-a).

Cluster 1 represents climatic series leading to low
yields from all locations. Clusters 3, 4, 6, 8 correpond rather to wheater types from the
north of France leading to high yields (3, 4, 6) or not (8).
 Clusters 2, 4, 6, 9 can be characterized by a cold weather and high yields
 but there are difficult to distinguish from each other; there is indeed an
 important mixture of clusters in node 6 in Figure
 \ref{fig:classif_caracterisation}-c and one can make the same
observation when studying the projection in Figure \ref{fig:classif_caracterisation}-a.

While a simple characterisation of clusters can be done, there are still
differences between them that we do not achieve to characterize,
which motivates the approach of using a distance between time series.
Especially, there is a known high impact of rain episodes and their
localization in time, however, the temporal mean aggregation of rain
is not retained when building these decision trees.

\subsection{Phenotype optimization}

\subsubsection{Algorithm performance}
\begin{figure}
\begin{center}
\includegraphics[width=.8\textwidth]{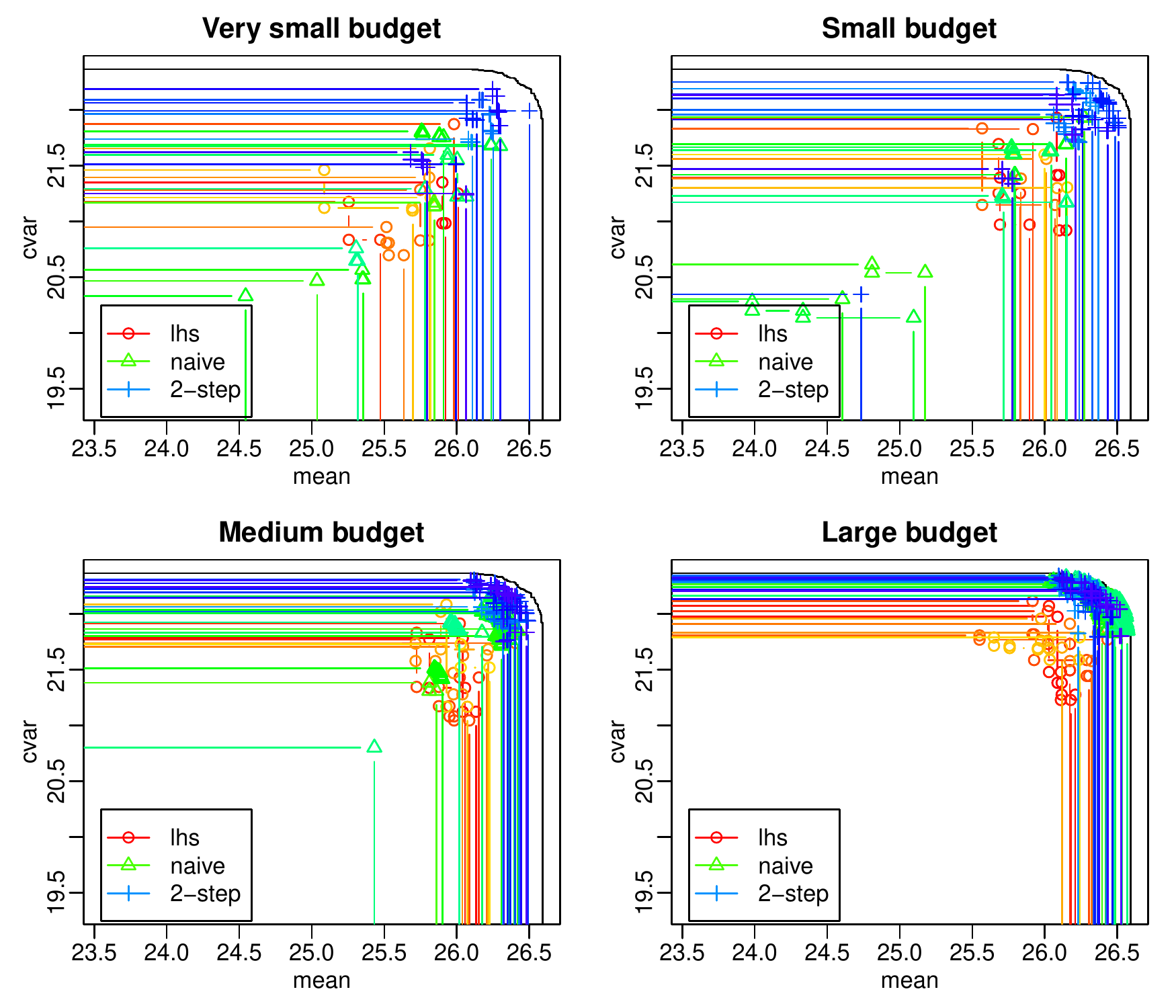}
\caption{\label{fig:allPareto} Pareto fronts obtained with the different methods for the four budgets considered.}
\end{center}
\end{figure}

\begin{figure}
\begin{center}
\includegraphics[width=\textwidth]{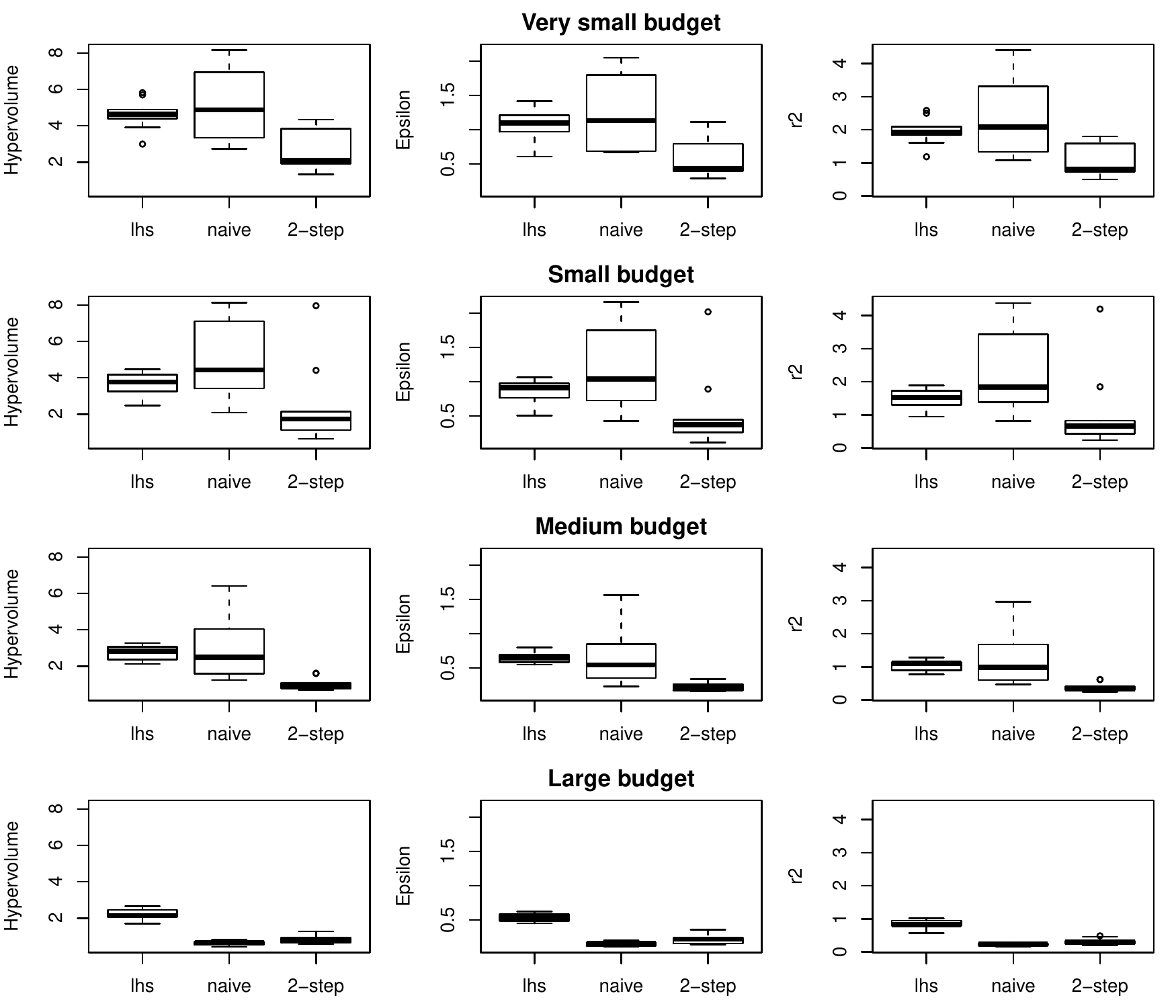}
\caption{\label{fig:allBoxPlots} Performance indices of the different methods for the four budgets considered.}
\end{center}
\end{figure}

Next, we compare the performances of the three approaches. As measuring performance
is non-trivial in multi-criteria optimization, we use three indicators: hypervolume,
epsilon and $R_2$ indicators \citep[as recommended in ][]{zitzler2003performance,hansen1998evaluating},
all available in the \texttt{R} package \texttt{emoa} \citep{mersmann2012emoa}. They provide
different measures of distance to the exact Pareto set and coverage of the objective space.
In short, the hypervolume indicator is a measure of the volume contained between the Pareto
front and a reference point (here, the worst value of each objective).
The epsilon indicator is a maximin distance between two Pareto fronts (here, we use
the exact Pareto front as reference), while the $R_2$ indicator can be seen as an
average distance. Figure \ref{fig:allPareto} shows all the Pareto fronts (of the different runs and methods)
for the different budgets, and Figure \ref{fig:allBoxPlots} shows the corresponding
performance indicators in the form of boxplots.

For the very small budget, we see that no method succeeds at finding the exact Pareto front.
Besides, most of the Pareto fronts consists of a single point.
However, the two-step approach still largely outperforms random search, while a naive use of MOPSO-CD performs worst,
as it requires a certain number of iterations to find a descent direction.

For the small and medium budgets, the two-step approach consistently finds a good approximation of the Pareto front
(with the exception of two outliers with the small budget). For the three indicators, it clearly outperforms the
other approaches.

For the large budget, we see that the regular MOPSO-CD performs slightly better, which is expected.
Indeed, as soon as there is no necessity of parcimony, using approximate objectives instead of actual ones
tends to slower, rather than accelerate, convergence.

\subsubsection{Results analysis}
\begin{figure}
\begin{center}
\includegraphics[width=.8\textwidth]{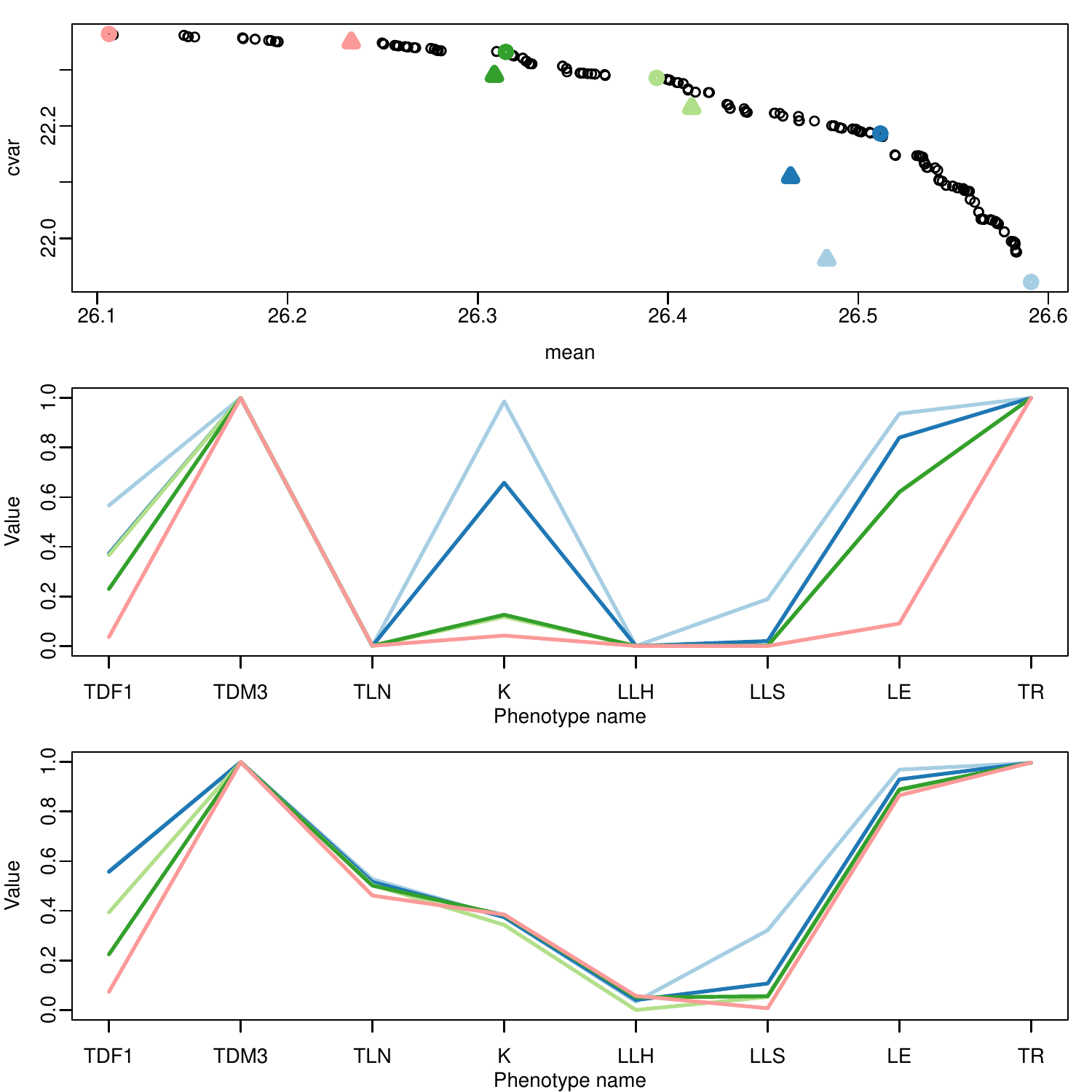}
  \caption{\label{fig:ParetoSet} Top: exact Pareto front. The bold circles correspond to a subset of five optimal phenotypes; the triangles correspond to five phenotypes returned by the two-step method. Middle: optimal phenotype values (one curve corresponds to one phenotype). Bottom: phenotype values obtained with the two-step method.}
\end{center}
\end{figure}

Finally, we characterize the results on the phenotype space. We compare here the exact Pareto set with one run
of the two-step method; we chose the run on the medium budget with the median performance. For readability, we
only consider a subset of the Pareto set of size five, equally spaced along the Pareto front. The Pareto fronts and sets
are represented in Figure \ref{fig:ParetoSet}.

We can see first that considering both the expectation and CVaR for optimization leads to a large variety of 
optimal phenotypes.
Looking back at the plant characteristics corresponding to those solutions, the optimum value
for five traits had little variability, meaning that those traits were important plant characteristics 
for crop performance in the tested environments.
Those five traits depicted plants adapted to water deficit: a late maturity (TDM3), a low leaf number (TLN),
largest leaves at the bottom of the plant (LLH), a small plant area (LLS), and a conservative strategy for 
stomatal conductance regulation (TR).
The three other traits (TDF1, K, LE) displayed variability in optimal values,
which was identified as the basis of the performance/stability trade-off (expectation/CVaR).
Here, the traits vary monotonically along the Pareto front.

Four distinct plant types could be identified in the phenotype space. For example, the \emph{red} plant type 
had an early flowering (TDF1), a low light extinction efficiency (K) and a low plant leaf area (LLS); 
those characterictics correspond to a conservative resource management strategy.
In an opposite manner, the \emph{light-blue} type displays a late flowering, a high efficiency to intercept 
light and a larger plant leaf area, characteristics usually associated with a productive but risky crop type 
when facing strong water deficit \citep{Connor1997}. The strategy associated with plant types identified 
from the phenotype space matched their position in the Pareto front,
i.e the \emph{light-blue} plant type was more performant but less stable than the \emph{red} one.

The Pareto set obtained with the two-step method reproduces part of these features: the fixed traits are 
similar (except TLN, which is fixed to approximately 0.5 instead of 0, (but this parameter is known to 
have little impact on the yield, see \citet{casadebaig.11}) and the variation of TDF1 and LLS is well-captured.
However, on this run the method failed at finding the variation of the K and LE traits: this probably explains
why the largest mean values (left of the Pareto front) are missed.

Overall, the two-step method allowed to identify the few key traits were responsible for the cultivar global 
adaptation capacity whereas secondary traits supported alternative resource use strategies underlying 
the yield expectation/stability tradeoff.

\section{Summary and perspectives}\label{sec:conclusion}
In this article, we proposed an algorithm for phenotype optimization under climatic uncertainties.
Our approach does not require any \textit{a priori} knowledge on the system besides parameter bounds,
hence is usable with any simulator depending on similar climatic data. Using subset selection for
the climates allowed us to reduce substantially the computational time without adding implementation issues.
If bias correction seems inevitable during optimization, we showed that a two-step strategy was sufficient
to achieve convergence: this point is critical as it allows our approach to be combined with any black-box
multi-objective solver.

Nevertheless, we see many opportunities for further improvements. First, the distance used here between climate series
does not account for the fact that agronomical systems are mostly sensitive to a few critical periods
(e.g., during flowering, grain filling). Weighting the DTW distance using expert knowledge or the results of a sensitivity
analysis may greatly improve the classification of the climates with respect to their impact on the model.

Second, the reconstruction step may benefit from additionnal study, in particular the effect of the subset size,
which has been fixed to 10 in our study for practical reasons but could be chosen using preliminary experiments for instance.
Another interesting topic would be to target the reconstruction to improve the quality of the objectives.
Indeed, the proposed approach aims a reconstructing the entire output distribution, while it is only important to obtain
good estimates of the expectation and the CVaR.

Third, a popular strategy to reduce the computational costs is to combine optimization with the use of surrogate modelling
\citep[see for instance][for recent examples]{di2009efficient,tsoukalas2014multiobjective}.
Our approach straightforwardly extends to such approaches, and would result in very parcimonious
algorithms that may be beneficial for expensive simulations.

Finally, we have chosen here to use a two-step strategy to allow the use of ``off-the-shelf'' optimization solvers.
Interlinking optimization and learning may improve substantially the efficiency of the method, although requiring
the development of an \textit{ad hoc} algorithm.

\section*{Appendix: clustering algorithm}\label{sec:appendix}
This section details the clustering algorithm used and the rule to chose the representative element of each class. 
The key of this particular approach is that, contrarily to a standard k-means algorithm, we cannot compute explicitely 
a central element (i.e., a ``virtual'' climatic series).

\begin{algorithm}[H]
 \caption{Clustering algorithm}
 \label{al:kmeans}
 \begin{algorithmic}
  \State Initialize $\boldsymbol{\beta}$ in $\Rset^{N \times K}$ randomly such that $\beta_{ij} \geq 0$, $\forall i,j $ and $\sum_{i=1}^N \beta_{ij} = 1$, $\forall j$. 
  Each line $\boldsymbol{\beta}_k$ is the dissimilarity of the centroid $\tilde \omega_k$ to the climates.
  \For {t= 1, \ldots, T}
  \State Pick $i$ randomly in $1, \ldots, N$ (one climate selected randomly)
  \State \textit{Assignment step} Find $j$ (center closest to $c_i$) such that 
  $$ j = \arg \min_{k = 1, \ldots, K} \left(\boldsymbol{\beta}_k^T \boldsymbol{\Delta}_i \right) - \frac{1}{2} \boldsymbol{\beta}_k \boldsymbol{\Delta} \boldsymbol{\beta}_k^T,$$
  with $\boldsymbol{\Delta}_i$ the i-th line of $\boldsymbol{\Delta}$.
  \State \textit{Representation step} (update center) $$\boldsymbol{\beta}_j \gets \boldsymbol{\beta}_j + r(t) \times (\one_j - \boldsymbol{\beta}_j),$$
  where $\one_j$ is a vector of zeros except its j-th value equal to one and $r(t) = \frac{\epsilon_0}{1 + c_0 \frac{t}{K}}$.
  \EndFor
 \end{algorithmic}
\end{algorithm}

Once $\boldsymbol{\beta}$ has converged, each climate $c_i$ is attributed to the class $j$, using:
\begin{equation*}
 j = \arg \min_{k = 1, \ldots, K} \left(\boldsymbol{\beta}_k \boldsymbol{\Delta}_i \right) - \frac{1}{2} \boldsymbol{\beta}_k \boldsymbol{\Delta} \boldsymbol{\beta}_k^T.
\end{equation*}

For each class $k$, a representative element $\omega^k$ is chosen.
We choose here the most central element in terms of dissimilarity. Let $\boldsymbol{\Delta}^k$ 
be the submatrix of $\boldsymbol{\Delta}$ corresponding to the elements of $\class^k$. We choose:
\begin{equation*}
 \omega^k = c_I^k \quad \text{ with } I=\arg \min_{1 \leq i \leq N^k} \sum_{j=1}^{N^k} \delta_{ij}^k
\end{equation*}


\end{document}